\documentstyle{article}
\newcommand{\be}{\begin{equation}}
\newcommand{\ee}{\end{equation}}
\newcommand{\ba}{\begin{array}}
\newcommand{\ea}{\end{array}}
\newcommand{\beqn}{\begin{eqnarray}}
\newcommand{\eeqn}{\end{eqnarray}}
\newcommand{\beqna}{\begin{eqnarray*}}
\newcommand{\eeqna}{\end{eqnarray*}}

\newcommand{\zero}{\setcounter{equation}{0}}

\font\tenbl=msbm10
\font\sevenbl=msbm7
\font\fivebl=msbm5
\newfam\blfam

\textfont\blfam=\tenbl
\scriptfont\blfam=\sevenbl
\scriptscriptfont\blfam=\fivebl

\font\symb=msam7
\def\znakr{\raise1.5pt\hbox{\symb\char66\kern-2pt\char74}}
\def\znakl{\raise1.5pt\hbox{\symb\char73\kern-2pt\char67}}

\begin{document}
\title{The $\kappa -$Weyl group and its algebra}
\author{Piotr Kosi\'nski\thanks{Supported by KBN grant 2P30221706p02} \\
Department of Theoretical Physics \\
University of \L \'od\'z \\
ul.Pomorska 149/159, 90-236 \L \'od\'z , Poland
\and
Pawel Ma\'slanka${}^{*}$\\
Department of Functional Analysis \\
University of \L \'od\'z \\
ul.Banacha 22, 90-238 \L \'od\'z , Poland}
\date{}
\maketitle
\begin{abstract}
The $\kappa -$Poincare group and its algebra in an arbitrary basis
are constructed. The $\kappa -$de\-formation of the Weyl group and its
algebra in any dimentions and in the reference frame in which $g_{00}=0$
are discussed.
\end{abstract}
\setcounter{section}{0}

\section{Introduction}
\zero

It is our great pleasure to contribute to this volume.

In last years we had an opportunity to collaborate with Prof. Jerzy
Lukierski. Our common topic has been the deformed symmetries of
space-time, mainly the so called $\kappa -$Poincare algebra invented
by Lukierski, Nowicki and Ruegg${}^1$. Apart from investigating the formal
properties of $\kappa -$Poincare algebra and looking for its possible
physical applications one of the main ideas of Prof. Lukierski is to extend the
notion of $\kappa$-deformation to larger groups of space-time symmetries.
This idea resulted in series of papers${}^{2-6}$ devoted to the
$\kappa -$deformation of SUSY extensions of the Poincare symmetry.
The next step to be done is to look for $\kappa -$deformed conformal
group / algebra. This problem has not been fully solved yet but some
perliminary steps were already undertaken${}^{7,8}$.

Inspired by these papers and numerous discussions with Prof. Lukierski we
attempt here to make a small step toward the solution of this problem.

Classically, the conformal group in four dimensions is nothing but
$SO(4,2)$. However, the standard (martix) parametrization of $SO(4,2)$
is not used, when $SO(4,2)$ is viewed as conformal group. On the contrary,
the conformal group is obtained from the action of $SO(4,2)$ on light
cone in sixdimensional space-time. But the light-cone coordinates are
related in rather complicated way to Minkowski coordinates in four
dimensions. Consequently, the standard parametrization of  $SO(4,2)$ is
related to the ``conformal" one by a complicated (even somewhere singular)
change of group parameters. This poses no problem on the ``classical" level.
However, if we are passing to the ``quantum" (i.e. deformed) case we are
faced with typical ordering problems of quantization procedure. This gives
some flavor of difficulties one meets trying to deform the conformal group.

In the recent paper${}^8$, Lukierski, Minnaert and Mozrzymas considered a new
class of classical $r -$matrices on conformal algebras in three and four
dimensions, which obey the classical Yang-Baxter equation and depend on
dimensionfull parameter. An important observation concerning the $d=4$
case was that the classical $r -$matrices obtained by them depend only on
genetators belonging to Poincare subalgebra of conformal algebra. Due to the
fact that they obey the classical Yang-Baxter equation (and not modified one)
they provide $r -$matrices for any algebra, containing Poincare algebra as
subalgebra. One of the $r -$matrices considered in${}^8$ leads to the so called
null-plane deformation of Poincare algebra found, by different methods, in
Ref.${}^9$. This deformation is similar to the standard $\kappa -$deformation.
The only diffrence is in the choice of undeformed subalgebra which is the
stability algebra of light-like fourvector instead of time-like one.
However, this difference is significant: in the standard case the Schouten
bracket is ad-invariant but does not vanish. Therefore, the relevant
$r -$matrix does not provide automatically the $r -$matrix for any extension
of Poincare algebra. Actually, the invariance
is broken already after adjoining the dilatation generator $D$.
Our aim here is to put the results of Ref.${}^8$ in more general setting.
In sec.\ref{r2} we review the properties of Poincare group for arbitrary
chosen mertic and discuss the Poisson structure on it. In sec.\ref{r3} the
quantization of this classical structure is performed. The bicrossproduct
form of resulting quantum group allows us to find, by duality, the relevant
algebra. The Weyl group and its algebra are constructed in sec.\ref{r4}.
Finally, sec.\ref{r5} is devoted to some conclusions.

\section{Relativity theory in on arbitrary coordinate system. \label{r2}}
\zero

Let us consider the $n -$dimensional lineral metric space $M$ with metric
tensor $g_{\mu \nu \ (\mu ,\nu =0,1,...,n-1)}$ given by an arbitrary
nondegenerate symmetric $n\times n$ matrix (not necessery diagonal).

Poincare group ${\cal P}$ is the group of inhomogenous transformations of the
space $M$:
$$
x^{\prime \mu }= \Lambda ^\mu _{\ \nu } x^\nu +a^\mu
$$
where the matrices $\Lambda ^\mu _{\ \nu }$ (Lorentz group) satisfy the
condition:
$$
g_{\mu \nu }= \Lambda ^\alpha _{\ \mu } \Lambda ^\beta _{\ \nu }
g_{\alpha \beta } .
$$

It is easy to see that the Poincare algebra $\tilde {{\cal P}}$ reads:
\beqn
\lbrack  P_\mu  , P_\nu \rbrack  & = & 0 \nonumber \\
\lbrack  M_{\mu \nu } , P_\lambda \rbrack  & = & i (g_{\nu \lambda }
P_\mu - g_{\mu \lambda } P_\nu ) \nonumber \\
\lbrack  M_{\mu \nu } , M_{\lambda \sigma } \rbrack  & = & i(g_{\mu \sigma }
M_{\nu \lambda }-g_{\nu \sigma }M_{\mu \lambda }+g_{\nu \lambda }
M_{\mu \sigma }-g_{\mu \lambda }M_{\nu \sigma }) \nonumber
\eeqn
where
$$
(M_{\alpha \beta })^\mu _{\ \nu } = i (\delta ^\mu _{\ \alpha }
g_{\nu \beta }-\delta ^\mu _{\ \beta }g_{\nu \alpha }) .
$$

Now consider $r \in \bigwedge ^2 \tilde {{\cal P}}$ given as follows,
${}^{10}$:
\be
r = {i\over \kappa } M_{0\nu } \wedge P^\nu  =
r^{\mu \nu ,\alpha }M_{\mu \nu }\wedge P_\alpha \label{rmat}
\ee
where
$$
r^{\mu \nu ,\alpha }={i\over 2\kappa }(\delta ^\mu _{\ 0}
g^{\nu \alpha }-\delta ^\nu _{\ 0}g^{\mu \alpha })
$$
and $\kappa $ is a real deformation parameter.

A calculation of Schouten bracket of $r$ with itself yields
\be
\lbrack r,r \rbrack = {ig_{00}\over \kappa ^2}M_{\alpha \beta }
\wedge P^\alpha \wedge P^\beta . \label{komr}
\ee

It is not dificult to see that $\lbrack r,r \rbrack $ is invariant,
hence $r$ is defines a structure of a Poisson Lie group on ${\cal P}$,
by the formula:
\be
\{ f,g \} = 2r^{\alpha \beta }( X^R_{\ \alpha }fX^R_{\ \beta }g-
X^L_{\ \alpha }fX^L_{\ \beta }g) \label{npois}
\ee
where $X^R_{\ \alpha }, X^L_{\ \beta }$ are the right- and left-invariant
vector fields.

It is easy to find the following expressions for the invariant
vector fields:
\beqn
X^{\ \alpha \beta }_L &=& \Lambda ^{\mu \alpha }{\partial \over
\partial \Lambda ^\mu _{\ \beta }} - \Lambda ^{\mu \beta }
{\partial \over \partial \Lambda ^\mu _{\ \alpha }} \nonumber \\
X^{\ \alpha }_L &=& \Lambda ^{\mu \alpha }{\partial \over
\partial a^\mu }\nonumber \\
X^{\ \alpha \beta }_R &=& \Lambda ^\beta _{\ \nu }
{\partial \over \partial \Lambda _{\alpha \nu }}-
\Lambda ^\alpha _{\ \nu }{\partial \over
\partial \Lambda _{\beta \nu }}+
a^\beta {\partial \over \partial a_\alpha }-
a^\alpha {\partial \over \partial a_\beta }\nonumber \\
X^{\ \alpha }_R &=& {\partial \over \partial a_\alpha } .\nonumber
\eeqn

In the Lie algebra basis corresponding to the above vector fields
we have the following relation between the generators of the Poincare
algebra and the invariant vector fields:
\beqn
M^{\alpha \beta } &=& i X^{\alpha \beta } \nonumber \\
P^\alpha &=& X^\alpha . \nonumber
\eeqn

This enables us to calculate the Poisson brackets of the coordinate
functions ${\cal P}$:
\beqn
\{ \Lambda ^\alpha _{\ \beta },a^\varrho  \}&=& - {1\over \kappa }
((\Lambda ^\alpha _{\ 0}-\delta ^\alpha _{\ 0})
\Lambda ^\varrho _{\ \beta }+(\Lambda _{0\beta }-
g_{0\beta })g^{\alpha \varrho }) \nonumber \\
\{ a^\varrho ,a^\sigma  \} &=& {1\over \kappa }
(\delta ^\varrho _{\ 0}a^\sigma  -
\delta ^\sigma _{\ 0}a^\varrho ) \nonumber \\
\{ \Lambda ^\alpha _{\ \beta }, \Lambda ^\mu _{\ \nu } \} &=& 0 .\label{wzor1}
\eeqn
\section{The $\kappa -$Poincare group and $\kappa -$Poincare algebra
in an arbitrary basis. \label{r3}}
\zero

If we perform the standard quantizations of the Poisson brackets
of the coordinate functions on ${\cal P}$ by replacing $\{ \ ,\ \} \to
{1\over i}\lbrack \ ,\ \rbrack $ one obtains the following set of
commutation relations:
\beqn
\lbrack  \Lambda ^\alpha _{\ \beta },a^\varrho  \rbrack &=& - {i\over \kappa }
((\Lambda ^\alpha _{\ 0}-\delta ^\alpha _{\ 0})
\Lambda ^\varrho _{\ \beta }+(\Lambda _{0\beta }-
g_{0\beta })g^{\alpha \varrho }) \nonumber \\
\lbrack  a^\varrho ,a^\sigma  \rbrack  &=& {i\over \kappa }
(\delta ^\varrho _{\ 0}a^\sigma  -
\delta ^\sigma _{\ 0}a^\varrho ) \nonumber \\
\lbrack  \Lambda ^\alpha _{\ \beta }, \Lambda ^\mu _{\ \nu } \rbrack
&=& 0 . \label{komut}
\eeqn

This standard quantization procedure is unambiguous: there is no ordering
ambiguity when quantizing the right-hand side of Eg.(\ref{wzor1}) due
to the commutativity of $\Lambda $'s.

Since the composition law is compatible with Poisson brackets, the above
commutation rules are compatible with the following coproduct:
\beqn
\Delta \Lambda ^\mu _{\ \nu } &=& \Lambda ^\mu _{\ \alpha } \otimes
\Lambda ^\alpha _{\ \nu }\nonumber \\
\Delta (a^\mu ) &=& \Lambda ^\mu _{\ \nu }\otimes a^\nu  +
a^\mu \otimes I . \label{delta}
\eeqn
The antipode and the counit are given by:
\beqn
S(\Lambda ^\mu _{\ \nu }) &=& \Lambda ^{\ \mu }_\nu \nonumber \\
S(a^\mu )&=& - \Lambda ^{\ \mu }_\nu  a^\nu \nonumber \\
\varepsilon (\Lambda ^\mu _{\ \nu })&=& \delta ^\mu _{\ \nu }\nonumber \\
\varepsilon (a^\mu )&=& 0 . \label{siep}
\eeqn

If we define the $* -$operation in such away that
$\Lambda ^\mu _{\ \nu }$ and $a^\mu $ are  selfadjoint elements, we conclude
that the relations Eq.(\ref{komut}), (\ref{delta}), (\ref{siep}) define a
Hopf $*-$algebra --- the $\kappa -$Poincare group ${\cal P}_\kappa $.
It follows from the Eq.(\ref{komut}), (\ref{delta}), (\ref{siep})
that the form of the $\kappa -$Poincare group does not depend on the choice
of the matric tensor $g_{\mu \nu }$. The differences between the various
$\kappa -$Poincare group are related to the fact that $\Lambda _{0\beta }$
appearing in the first commutation relation of Eq.(\ref{komut}) are not
the independ variables but are the linear combinations of the independent
ones: $\Lambda _{0\beta }= g_{0\mu }\Lambda ^\mu _{\ \beta }$.

It should be stressed, that the $\kappa -$Poincare group can be defined
as a right-left bicrosspro\-duct${}^{11,12}$:
$$
{\cal P}_\kappa = T^* \znakr C(L) .
$$
To see this it is sufficient to define the structure maps:
\beqn
\beta (a^\mu ) &=& \Lambda ^\mu _{\ \nu } \otimes x^\nu  \nonumber \\
\Lambda ^\mu _{\ \nu } \lhd x^\varrho  &=&
- {i\over \kappa } (( \Lambda ^\mu _{\ 0} - \delta ^\mu _{\ 0})
\Lambda ^\varrho _{\ \nu } + (\Lambda _{0\nu }- g_{0\nu })
g^{\mu \varrho }) . \nonumber
\eeqn
Moreover, while $C(L)$ is the standard algebra of function defined over
Lorentz group, $T^*$ is defined by the following relations:
\beqn
\lbrack a^\mu , a^\nu   \rbrack  &=& {i\over \kappa }
(\delta ^\mu _{\ 0}a^\nu  - \delta ^\nu _{\ 0 }a^\mu ) \nonumber \\
\Delta (a^\mu ) &=& a^\mu  \otimes I + I\otimes a^\mu \nonumber \\
S(a^\mu ) &=& -a^\mu \nonumber \\
\varepsilon (a^\mu ) &=& 0 . \nonumber
\eeqn
The bicrossproduct structure of the $\kappa -$Poincare group allows us to
define the dual object, the $\kappa -$Poincare algebra
$\tilde {{\cal P}}_\kappa $
as a left-right bicrossproduct:
$$
\tilde {{\cal P}}_\kappa  = T \znakl U(\tilde L)
$$
where $T$ is dual to $T^*$ and $U(\tilde L)$ is the universal
enveloping algebra of the Lorentz algebra.
The duality $T^* \iff T$ is defined by:
$$
<a^\mu , P_\nu > = i\delta ^\mu _{\ \nu }.
$$
The duality between the Lorentz group and algebra is defined in the
standard way:
$$
<\Lambda ^\mu _{\ \nu },M^{\alpha \beta }>
= i(g^{\alpha \mu }\delta ^\beta _{\ \nu }-g^{\beta \mu }
\delta ^\alpha _{\ \nu }) .
$$
The structure maps are defined by the following duality relations:
\beqn
<t,M^{\alpha \beta }\rhd P_\gamma > &=& <\beta (t),M^{\alpha \beta }
\otimes P_\gamma >\nonumber \\
<\Lambda \lhd t,M^{\alpha \beta }> &=& <\Lambda \otimes t,
\delta (M^{\alpha \beta })>  \nonumber
\eeqn
here t is arbitrary product of a's while $\Lambda $ is an arbitrary product
of $\Lambda $'s.

Finally using the method described in${}^{13}$ we arrive at the following
explicite form of the  $\kappa -$Poin\-care algebra:

\noindent a) the commutation rules:
\beqn
\lbrack  M^{ij},P_0 \rbrack  &=& 0 \nonumber \\
\lbrack M^{ij},P_k  \rbrack  &=& i\kappa (\delta ^j_{\ k}
g^{0i}-\delta ^i_{\ k}g^{0j})(1-e^{-{P_0\over \kappa }})
+ i(\delta ^j_{\ k}g^{is}-\delta ^i_{\ k}g^{js})P_s\nonumber \\
\lbrack M^{i0},P_0  \rbrack  &=& i\kappa g^{i0}
(1-e^{-{p_0\over \kappa }})+ig^{ik}P_k \nonumber \\
\lbrack M^{i0},P_k  \rbrack  &=& -i{\kappa \over 2}g^{00}
\delta ^i_{\ k}(1-e^{-2{P_0\over \kappa }})-
i\delta ^i_{\ k}g^{0s}P_s e^{-{P_0\over \kappa }}+ \nonumber \\
 &\ & +ig^{0i}P_k (e^{-{P_0\over \kappa }}-1)+
{i\over 2\kappa }\delta ^i_{\ k}g^{rs}P_r P_s -
{i\over \kappa }g^{is} P_s P_k \nonumber \\
\lbrack P_\mu ,P_\nu   \rbrack  &=& 0 \nonumber \\
\lbrack M^{\mu \nu },M^{\lambda \sigma }  \rbrack
&=& i(g^{\mu \sigma }M^{\nu \lambda }-g^{\nu \sigma }M^{\mu \lambda }
+g^{\nu \lambda }M^{\mu \sigma }-g^{\mu \lambda }M^{\nu \sigma }) \nonumber
\eeqn
b) the coproducts:
\beqn
\Delta P_0 &=& I\otimes P_0 + P_0 \otimes I \nonumber \\
\Delta P_k &=& P_k \otimes e^{-{P_0\over \kappa }}+
I\otimes P_k \nonumber \\
\Delta M^{ij} &=& M^{ij}\otimes I + I\otimes M^{ij} \nonumber \\
\Delta M^{i0} &=& I\otimes M^{i0} + M^{i0}\otimes
e^{-{P_0\over \kappa }}- {1\over \kappa }M^{ij}\otimes P_j \nonumber
\eeqn
where $i,j,k = 1,2,3,...,n-1$.

Let us note that the $\kappa -$Poincare algebra and group in any
dimensions and for the diagonal metric tensor were obtained in${\ }^{14}$
and${\ }^{15}$. However the duality between the  $\kappa -$Poincare
algebra and group was not discussed in these papers.

In the end of this section let us remark that the classical $r -$matrix
Eq.(\ref{rmat}), $r=M_{0i}\wedge P^i$, does not modify the coproducts
for the generators $M^{ij}$, forming undeformed Lie subalgebra
as well as the component $P_0$ of the fourmomentum. The algebra with
the generators $M^{ij}$, $P_0$ describes the classical subalgebra of our
$\kappa -$Poincare algebra.
\section{The $\kappa -$Weyl group and algebra \label{r4}}
\zero

The classical Weyl group ${\cal W}$ consists of the triples $(a,\Lambda ,e^b)$,
where $a$ is a $n-$vector, $\Lambda $ is the matrix of the Lorentz
group in $n-$dimensions and $b\in R$, with the composition law:
$$
(a^\mu ,\Lambda ^\mu _{\ \nu },e^b)*(a^{\prime \nu },
\Lambda ^{\prime \nu }_{\ \alpha },e^{b^\prime }) =
(\Lambda ^\mu _{\ \nu }e^b a^{\prime \nu }+a^\mu ,
\Lambda ^\mu _{\ \nu }\Lambda ^{\prime \nu }_{\ \alpha },
e^b e^{b^\prime }) .
$$
Its Lie algebra, the Weyl algebra $\tilde {\cal W}$ reads:
\beqn
\lbrack P_\mu ,P_\nu   \rbrack  &=& 0 \nonumber \\
\lbrack M_{\mu \nu },P_\lambda   \rbrack  &=&
i(g_{\nu \lambda }P_\mu - g_{\mu \lambda }P_\nu ) \nonumber \\
\lbrack M_{\mu \nu },M_{\lambda \sigma }  \rbrack
&=& i(g_{\mu \sigma }M_{\nu \lambda }-g_{\nu \sigma }M_{\mu \lambda }
+g_{\nu \lambda }M_{\mu \sigma }-g_{\mu \lambda }
M_{\nu \sigma }) \nonumber \\
\lbrack M_{\mu \nu },D  \rbrack  &=& 0 \nonumber \\
\lbrack P_\mu ,D  \rbrack   &=& -iP_\mu  \nonumber
\eeqn
here $M_{\mu \nu }$, $P_\mu $ are the generators of the Poincare algebra
and $D$ is the dilatation generator.

We would like to obtain the $\kappa -$deformation of the Weyl group and its
Lie algebra. It is clear that the classical $r-$matrix for the Poincare
algebra satisfyning the classical Yang-Baxter equation (CYBE) is also
a classical $r-$matrix for the Weyl algebra. From the Eq.(\ref{komr})
follows that our classical $r-$matrix Eq.(\ref{rmat}) satisfies the CYBE
iff $g_{00}=0$. This means that our $r-$matrix defines a structure of
a Poisson Lie group on the Weyl group only in the basis in which the
metric tensor takes such a from that $g_{00}=0$. We shall consider only
these types of metrices.

In order to obtain the $\kappa -$deformation of the Weyl group we firstly
find the invariant fields:
\beqn
X^{\ \alpha \beta }_L &=& \Lambda ^{\mu \alpha }{\partial \over
\partial \Lambda ^\mu _{\ \beta }} - \Lambda ^{\mu \beta }
{\partial \over \partial \Lambda ^\mu _{\ \alpha }} \nonumber \\
X^{\ \alpha }_L &=& e^b \Lambda ^{\mu \alpha }{\partial \over
\partial a^\mu }\nonumber \\
X_L &=& e^b {\partial \over \partial e^b}\nonumber \\
X^{\ \alpha \beta }_R &=& \Lambda ^\beta _{\ \nu }
{\partial \over \partial \Lambda _{\alpha \nu }}-
\Lambda ^\alpha _{\ \nu }{\partial \over
\partial \Lambda _{\beta \nu }}+
a^\beta {\partial \over \partial a_\alpha }-
a^\alpha {\partial \over \partial a_\beta }\nonumber \\
X^{\ \alpha }_R &=& {\partial \over \partial a_\alpha } \nonumber \\
X_R &=& a^\mu {\partial \over \partial a^\mu }+
e^b {\partial \over \partial e^b} . \nonumber
\eeqn

Then, using the Eq.(\ref{npois}) we calculate the Poisson brackets of the
coordinate functions on ${\cal W}$ and perform the standard quantizations of
the
Poisson brackets, by replacing $\{ \ , \ \} \to {1\over i}
\lbrack \ ,\ \rbrack $.

Finally we obtain the following set of commutation relations:
\beqn
\lbrack  \Lambda ^\alpha _{\ \beta },a^\varrho  \rbrack &=& - {i\over \kappa }
((e^b \Lambda ^\alpha _{\ 0}-\delta ^\alpha _{\ 0})
\Lambda ^\varrho _{\ \beta }+(\Lambda _{0\beta }-
e^b g_{0\beta })g^{\alpha \varrho }) \nonumber \\
\lbrack  a^\varrho ,a^\sigma  \rbrack  &=& {i\over \kappa }
(\delta ^\varrho _{\ 0}a^\sigma  -
\delta ^\sigma _{\ 0}a^\varrho ) \nonumber \\
\lbrack  \Lambda ^\alpha _{\ \beta }, \Lambda ^\mu _{\ \nu } \rbrack
&=& 0  \nonumber \\
\lbrack \Lambda ^\mu _{\ \nu },b \rbrack  &=& 0 \nonumber \\
\lbrack a^\mu ,b \rbrack  &=& 0 . \nonumber
\eeqn
This standard quantization procedure is unambiguonus: there is no
ordering ambiguity. Since the composition law is compatible with
the Poisson brackets, the above commutation rules are compatible
with the following coproduct:
\beqn
\Delta \Lambda ^\mu _{\ \nu } &=& \Lambda ^\mu _{\ \alpha } \otimes
\Lambda ^\alpha _{\ \nu }\nonumber \\
\Delta a^\mu  &=& e^b \Lambda ^\mu _{\ \nu }\otimes a^\nu  +
a^\mu \otimes I   \nonumber \\
\Delta b &=& b\otimes I + I\otimes b  . \label{delta1}
\eeqn
The antipode and the counit are given by:
\beqn
S(\Lambda ^\mu _{\ \nu }) &=& \Lambda ^{\ \mu }_\nu \nonumber \\
S(a^\mu )&=& - e^{-b} \Lambda ^{\ \mu }_\nu  a^\nu \nonumber \\
S(b) &=& -b \nonumber \\
\varepsilon (\Lambda ^\mu _{\ \nu })&=& \delta ^\mu _{\ \nu }\nonumber \\
\varepsilon (a^\mu )&=& 0  \nonumber \\
\varepsilon (b) &=& 0 . \label{siep1}
\eeqn
We conclude that the Eq.(\ref{delta1}), (\ref{siep1}) define
the Hopf algebra --- $\kappa -$Weyl group ${\cal W}_\kappa $.

If we forget for a moment about our general theory, we can check by
explicite calculation that our structure is self consistent
(Jacobi identities, the relations $\lbrack \Delta a,\Delta b \rbrack =
\Delta  \lbrack a,b  \rbrack $) iff $g_{00}=0$ or $b=0$.

For example:
$$
\lbrack \Lambda ^\alpha _{\ \beta },\lbrack a^\varrho ,a^\sigma
\rbrack \rbrack + \lbrack a^\sigma ,\lbrack \Lambda ^\alpha _{\ \beta },
a^\varrho \rbrack \rbrack + \lbrack a^\varrho ,\lbrack a^\sigma ,
\Lambda ^\alpha _{\ \beta }\rbrack \rbrack =
 {1\over \kappa ^2}g_{00}(1-e^{-2b})(g^{\alpha \sigma }
\Lambda ^\varrho _{\ \beta }- g^{\alpha \varrho }
\Lambda ^\sigma _{\ \beta }).
$$
It is easy to see that the ${\cal W}_\kappa $ has a right-left bicrossproduct
structure:
$$
{\cal W}_\kappa = T^* \znakr C(A)
$$
where $C(A)$ is the standard algebra of functions defined over group $A$.
The group $A$ consists of the pairs $(\Lambda ,e^b)$, where $\Lambda $
is a matrix of the Lorentz group and $b\in R$ and with the composition
law:
$$
(\Lambda ,e^b)*(\Lambda ^\prime ,e^{b^\prime }) = (\Lambda \Lambda ^\prime ,
e^b e^{b^\prime }) .
$$
$T^*$ is defined by the rerations:
\beqn
\lbrack  a^\varrho ,a^\sigma  \rbrack  &=& {i\over \kappa }
(\delta ^\varrho _{\ 0}a^\sigma  -
\delta ^\sigma _{\ 0}a^\varrho ) \nonumber \\
\Delta a^\mu  &=& a^\mu  \otimes I + I\otimes  a^\mu  \nonumber \\
S(a^\mu ) &=& -a^\mu  \nonumber \\
\varepsilon (a^\mu ) &=& 0 . \nonumber
\eeqn
To see this it is sufficient to define the structure maps as follows:
\beqn
\beta (a^\mu ) &=& e^b \Lambda ^\mu _{\ \alpha }\otimes a^\alpha \nonumber \\
\Lambda ^\mu _{\ \nu }\lhd a^\varrho  &=& -{i\over \kappa }
((e^b \Lambda ^\mu _{\ 0} - \delta ^\mu _{\ 0})
\Lambda ^\varrho _{\ \nu } + (\Lambda _{0\nu }-
e^b g_{0\nu })g^{\mu \varrho }) \nonumber \\
e^b \lhd a^\varrho &=& 0 . \nonumber
\eeqn
This right-left bicrossproduct structure of the $\kappa -$Weyl group
allows us to define the $\kappa -$Weyl algebra $\tilde {\cal W}_\kappa $
as a left-right bicrossproduct structure:
$$
\tilde {\cal W}_\kappa  = T \znakl U(\tilde A )
$$
where $T$ is dual to $T^*$ and $U(\tilde A )$ is the universal enveloping
algebra of the Lie algebra $\tilde A$ of the group $A$.

The duality $T^* \iff T$ is defined by:
$$
<a^\mu ,P_\nu > = i\delta ^\mu _{\ \nu } .
$$
The duality between the group $A$ and the algebra $\tilde A$ is defined in the
standard way:
\beqn
<\Lambda ^\mu _{\ \nu },M^{\alpha \beta }> &=& i(g^{\alpha \mu }
\delta ^\beta _{\ \nu }- g^{\beta \mu }\delta ^\alpha _{\ \nu }) \nonumber \\
<b,M^{\alpha \beta }> &=& 0 \nonumber \\
<\Lambda ^\mu _{\ \nu },D> &=& 0 \nonumber \\
<b^n ,D> &=& i\delta _{n,1} . \nonumber
\eeqn
The structure maps are defined by the following duality relations:
\beqn
<t,M^{\alpha \beta }\rhd P_\gamma > &=& <\beta (t),M^{\alpha \beta }
\otimes P_\gamma >\nonumber \\
<t,D\rhd P_\gamma > &=& <\beta (t),D\otimes P_\gamma > \nonumber \\
<\Gamma \lhd t,M^{\alpha \beta }> &=& <\Gamma \otimes t,
\delta (M^{\alpha \beta })> \nonumber \\
<\Gamma \lhd t,D> &=& <\Gamma \otimes t,\delta (D)> \nonumber
\eeqn
here t is arbitrary product of a's while $\Gamma $ is an arbitrary product
of $\Lambda $'s.

Finally using the method described in${}^{13}$ we arrive at the following
explicite form of the  $\kappa -$Weyl algebra:

\noindent a) the commutation rules:
\beqn
\lbrack  M^{ij},P_0 \rbrack  &=& 0 \nonumber \\
\lbrack M^{ij},P_k  \rbrack  &=& i\kappa (\delta ^j_{\ k}
g^{0i}-\delta ^i_{\ k}g^{0j})(1-e^{-{P_0\over \kappa }})
+ i(\delta ^j_{\ k}g^{is}-\delta ^i_{\ k}g^{js})P_s\nonumber \\
\lbrack M^{i0},P_0  \rbrack  &=& i\kappa g^{i0}
(1-e^{-{p_0\over \kappa }})+ig^{ik}P_k \nonumber \\
\lbrack M^{i0},P_k  \rbrack  &=& -i{\kappa \over 2}g^{00}
\delta ^i_{\ k}(1-e^{-2{P_0\over \kappa }})
-i\delta ^i_{\ k}g^{0s}P_s e^{-{P_0\over \kappa }}+ \nonumber \\
 &\ & +ig^{0i}P_k (e^{-{P_0\over \kappa }}-1)+
{i\over 2\kappa }\delta ^i_{\ k}g^{rs}P_r P_s -
{i\over \kappa }g^{is} P_s P_k \nonumber \\
\lbrack D,P_0  \rbrack  &=& i\kappa (1-e^{-{P_0\over \kappa }}) \nonumber \\
\lbrack D,P_i  \rbrack  &=& iP_i e^{-{P_0\over \kappa }}
+i{\kappa \over 2}g^{00}g_{i0}(1-e^{-{P_0\over \kappa }})^2
+ig_{0i}g^{0s}P_s (1-e^{-{P_0\over \kappa }})
+{i\over 2\kappa }g_{0i}g^{rs}P_r P_s \nonumber
\eeqn
b) the coproducts:
\beqn
\Delta D &=& D\otimes I + I\otimes D -
g_{0i}M^{i0}\otimes (1-e^{-{P_0\over \kappa }})-
{1\over \kappa }g_{0i}M^{ik}\otimes P_k \nonumber \\
\Delta P_0 &=& I\otimes P_0 + P_0 \otimes I \nonumber \\
\Delta P_k &=& P_k \otimes e^{-{P_0\over \kappa }}+
I\otimes P_k \nonumber \\
\Delta M^{ij} &=& M^{ij}\otimes I + I\otimes M^{ij} \nonumber \\
\Delta M^{i0} &=& I\otimes M^{i0} + M^{i0}\otimes
e^{-{P_0\over \kappa }}- {1\over \kappa }M^{ij}\otimes P_j . \nonumber
\eeqn

\section{Conclusions \label{r5}}
\zero

We have constructed the  $\kappa -$Poincare group resulting from Poincare
group formulated in an arbitrary basis. The quantization is unambigous
due to the absence of ordering problems. The resulting quantum group has
a bicrossproduct structure. Using this and the methods developed in${}^{13}$
we were able to construct the relevant $\kappa -$Poincare algebra.
The Schouten bracket of the classical $r -$matrix we have used appeared
to be proportional to the $g_{00}$ component of the metric tensor.
Therefore in the reference frame choosen in such a way that $g_{00}=0$
the relevant Poisson structure can be extended to any group containing
Poincare group as a subgroup. This was used to define the Poisson
structure on the Weyl group, which allowed us to construct
$\kappa -$deformation of this group. Again we obtained a bicrossproduct
structure which allowed us to construct the relevant $\kappa -$Weyl
algebra.

The above construction seems to us to be a proper introductory step toward
the definition $\kappa -$deformed conformal group. One can attempt to
quantize the Poisson structure on conformal group resulting from the same
$r -$martix we used in the case of the Weyl group hoping that the ordering
problems could be overcome in some way (as for example in the case of
$\kappa -$Poincare supergroup${}^{2-6}$). As a next step one try to construct
the relevant algebra. This might be more difficult as the bicrossproduct
structure is lacking in the case of conformal group.

An alternating way of attacking the problem would be to try to incorporate
on the quantum level the property of conformal group that it can be
obtained from Weyl group by adding (in a special way) the operation
of inversion.

\section{References}
\zero
\

1. J.Lukierski, A.Nowicki, H.Ruegg, {\it Phys.Lett.} {\bf B293},
(1993), 419.

2. P.Kosi\'nski, J.Lukierski, P.Ma\'slanka, J.Sobczyk,
{\it Journ. of Phys. }{\bf A27}, (1994), 6827.

3. P.Kosi\'nski, J.Lukierski, P.Ma\'slanka, J.Sobczyk,
{\it Journ. of Phys. }{\bf A28}, (1995), 2255.

4. P.Kosi\'nski, J.Lukierski, P.Ma\'slanka, J.Sobczyk,
preprint JFTUWr883/94, hep 94-12114.

5. P.Kosi\'nski, J.Lukierski, P.Ma\'slanka, J.Sobczyk,
preprint JFTUWr887/95, q-alg. 95-01010.

6. P.Kosi\'nski, J.Lukierski, P.Ma\'slanka, J.Sobczyk,
``Quantum deformation of the Poincare supergroup",
in:``Quantum Groups. Formalism and applications".
{\it The Procedings of XXX Karpacz Winter School of Theoretical
Physics}, eds. J.Lukierski, Z.Popowicz, J.Sobczyk, PWN 1995, p.353.

7. M.Klimek, J.Lukierski, {\it Acta Phys.Pol. }{\bf B26} (1995), 1209.

8. J.Lukierski, P.Minnaert, M.Mozrzymas, preprint CPTMB/PT/95-6.

9. A.Ballesteros, F.J.Herranz, M.A.del Olmo, M.Santander, preprint March 1995.

10. S.Zakrzewski, {\it Journ.of Phys.} {\bf A27} (1994), 2075.

11. P.Zaugg, preprint MIT-CTP, (1994).

12. S.Majid, H.Ruegg, {\it Phys.Lett.} {\bf B334}, (1994), 348.

13. P.Kosi\'nski, P.Ma\'slanka, preprint IMU\L  3/94.

14. J.Lukierski, H.Ruegg, {\it Phys.Lett.} {\bf B329}, (1994), 189.

15. P.Ma\'slanka, {\it Journ. of Phys. }{\bf A26}, (1993), L1251.
\end{document}